\documentclass[12pt]{amsart}
\usepackage{amssymb}

            \newcommand{\n}{{{\mathfrak n}}}
            \newcommand{\h}{{{\mathfrak h\,}}}
            \newcommand{\C}{{\mathbb C}}
            \newcommand{\s}{{\mathcal Sym}}
            \newcommand{\Sing}{{\rm Sing}}
            \newcommand{\0}{{\bf 0}}
            \newcommand{\uu}{{\bf u}}
            \newcommand{\vv}{{\bf v}}
            \newcommand{\w}{{\bf w}}
            \newcommand{\kk}{{\bf k}}
            \newcommand{\T}{{\bf t}}
            \newcommand{\bF}{{\bf F}}
\def\wh#1{\widehat{#1}}
\def\wt#1{\widetilde{#1}}
\newtheorem{thm}{Theorem}
\newtheorem {prop} {Proposition}

\newtheorem {cor} {Corollary}

\begin{document}
\title{On Bethe vectors in the  $sl_{N+1}$ Gaudin model}
\author[{}] {S.~Chmutov \and I.~Scherbak}
\begin{abstract}

The note deals with the Gaudin model associated with the tensor
product of  $n$  irreducible finite-dimensional $sl_{N+1}$-modules
marked by distinct complex numbers $z_1,\dots, z_n$.
The  Bethe Ansatz  is a method to construct common eigenvectors of
the Gaudin hamiltonians by means of chosen singular vectors in the
factors and  $z_j$'s. These vectors are called Bethe vectors.

The question if the Bethe vectors are non-zero vectors is open.
By the moment, the only way to verify that was
based on a relation to critical points of the master function of
the Gaudin model, and non-triviality of a Bethe vector was proved
only in the case when the corresponding critical point is
non-degenerate (\cite{ScV, MV1}). However  degenerate critical
points do appear  in the  Gaudin  model (\cite[Section12]{ReV}).

We believe that the Bethe vectors never vanish, and suggest an approach that
does not depend on non-degeneracy of the corresponding critical point.
The idea is for a Bethe vector to choose a suitable subspace in the weight 
space and to check that the projection of the Bethe vector to this subspace 
is non-zero. We apply this approach to  verify non-triviality of Bethe 
vectors in new examples.

\end{abstract}
\maketitle

\pagestyle{myheadings}
\markboth{S.~Chmutov and I.~Scherbak}
{On Bethe vectors in the  $sl_{N+1}$ Gaudin model}

\section{Introduction}\label{S1}
We study the Gaudin model of statistical mechanics associated with the Lie
algebra  $sl_{N+1}(\C)$. The {\it space of states} of the model is the tensor
product
$$
L=L_{\Lambda(1)}\otimes\dots\otimes  L_{\Lambda(n)}\, ,
$$
where $L_{\Lambda(j)}$  is a finite-dimensional irreducible  $sl_{N+1}$-module
with highest  weight $\Lambda(j)\,$, $1\leq j\leq n$\,.
For the standard notions of representation theory see~\cite{FuH}.

\medskip
In the Gaudin model, the modules  $L_{\Lambda(1)}, \dots, L_{\Lambda(n)}$  are
the spin spaces of $n$ particles located at distinct points
$z_1,\dots,z_n\in \C$.  Write $z=(z_1,\dots,z_n)$.
The {\it Gaudin hamiltonians} $H_1(z),\dots\,, H_n(z)$  are  mutually commuting
linear operators in $L$ which are defined as follows,
$$
H_j(z)=\sum_{i\neq j}\, \frac{C_{ij}}{z_j-z_i}\,,\ \  1\leq j\leq n\,,
$$
here $C_{ij}$ acts as the Casimir operator on factors $L_{\Lambda(i)}$ and
$L_{\Lambda(j)}$ of the tensor product and as the identity on all other factors.

\medskip
One of the main problems in the Gaudin model is simultaneous diagonalization
of the operators $H_1(z),\dots\,, H_n(z)$.
The Gaudin hamiltonians commute with the diagonal action of $sl_{N+1}$ in $L$,
therefore it is enough to find  common eigenvectors and the eigenvalues
in the subspace of  singular vectors of a given weight, for every weight.

\medskip
The algebraic Bethe Ansatz is a  method  to construct  such
vectors. The idea is to find some function  $\vv=\vv(\T)$ taking
values in the weight subspace ($\T$ is a multidimensional
auxiliary variable) and to determine a certain special value of
its argument, $\T^{(0)}$,  in such a way  that $\vv(\T^{(0)})$ is
a common eigenvector of the hamiltonians. The  equations on $\T$
which determine these  special values of the argument are called
{\it the Bethe equations}, and $\vv(\T^{(0)})$  is called {\it the
Bethe vector}. For the Gaudin model, the Bethe
equations and the function $\vv(\T)$  are written in \cite{FeFRe,
ReV, SV}. On Bethe vectors in  the Gaudin model see
also~\cite{G,FaT,Re}.

\medskip 
It was believed that for generic $z$ one can find an
eigenbasis in the subspace of singular vectors
consisting of Bethe vectors only. This is indeed the case
for the tensor products of $sl_2(\C)$-modules and for the tensor
products of several  copies of first and last fundamental
$sl_{N+1}$-modules ({\cite{ScV, MV1}). Recent results
of~\cite{MV2} show however that generically other eigenvectors have
to be present in eigenbases as well. These other vectors are in
some sense ``more degenerate'' than Bethe vectors, see~\cite{F1}
and especially~\cite[Section~5.5]{F2} discussing the ``degeneracies''.

\medskip 
In the Bethe Ansatz, two problems naturally arise: to find solutions
of the Bethe equations, and to check non-triviality of the
corresponding Bethe vectors. Both problems are open and seem to 
be difficult ones. On solutions to the Bethe equations in some
particular cases, see~\cite{V,ScV,MV1,Sc}.

\medskip
The present note is devoted to the question if Bethe vectors are 
non-zero vectors. By the moment, the only known way to verify that 
was extremely non-direct, via the so-called {\it master function}.
Namely, it appeared that the Bethe equations in the Gaudin model
form the critical point system of a certain function $S(\T;z)$,
here $\T$ is a multidimensional variable and $z=(z_1,\dots,z_n)$
is fixed, \cite{ReV}. Moreover, the norm of the Bethe vector
$\vv(\T^{(0)})$ with respect to some (degenerate) bilinear form on
the tensor product turned out to be the Hessian of $S(\T;z)$ at
the critical point $\T^{(0)}$; hence the  Bethe vectors
corresponding to non-degenerate critical points of the function
$S(\T;z)$ appeared to be non-zero vectors, \cite{V, MV1}. In this
way, the non-triviality of Bethe vectors has been checked for
generic $z$ in the case of tensor products of $sl_2(\C)$-modules
and in the case of tensor products of several copies of first and
last fundamental $sl_{N+1}$-modules, \cite{ScV, MV1}.

\medskip
It is known however, that for some values of  $z$  the master function does
have degenerate critical points; an example can be found in~\cite[Section~12]{ReV}.
Notice that in that example the corresponding Bethe vector is a non-zero vector
as well. We believe that the Bethe vectors are always non-trivial.

\medskip\noindent{\bf Conjecture.} \ \ {\it In the $sl_{N+1}(\C)$ Gaudin model,
every Bethe vector is non-zero, for any $z$. For some  values of $z$  the number
of Bethe vectors (i.e. of solutions to the Bethe equations, i.e. of critical
points of the master function) may decrease, but the Bethe vectors still are
non-zero.}

\medskip
We suggest a more direct approach that does not depend on
non-degeneracy of the corresponding critical point. 
The idea is  to project a Bethe vector to  a suitable subspace in 
the space  of singular vectors of
a given weight and to check that the projection is non-zero.

\medskip
We exploit this idea in some examples of tensor products of 
irreducible  finite-dimen-sional $sl_{N+1}$-modules. The case of
the tensor product of $n=2$ modules is special. First of all, in
this case {\it all} values of $z$ are  generic. Indeed, as it was
pointed out in~\cite[Section~5]{ReV}, for any fixed $z_1\neq z_2$ the
linear change of variables  $\uu=(\T-z_1)/(z_2-z_1)$ turns the
Bethe system on $\T$ with $z=(z_1,z_2)$ into the Bethe system on
$\uu$ with $z=(0,1)$. Next,  the  Gaudin hamiltonias
$H_1(0,1)=-H_2(0,1)$  are reduced to the Casimir operator, and
hence act in any irreducible submodule of the tensor product by
multiplication by a constant, i.e. {\it any} singular vector is
their common eigenvector. Finally, non-triviality of
a Bethe vector for  $n>2$  in many cases  can be deduced from
non-triviality of a certain set of Bethe vectors corresponding to
$n=2$ and $z=(0,1)$, by means of {\it iterated singular vectors}
introduced in \cite{ReV}; see~\cite{Sc} for a more detailed explanation.

\medskip
Let  $L=L_{\Lambda(1)}\otimes L_{\Lambda(0)}$ be the tensor product
of two  irreducible finite-dimensional $sl_{N+1}$-modules, where
$L_{\Lambda(1)}$ is marked by $z_1=1$ and  $L_{\Lambda(0)}$ by
$z_0=0$. Denote simple positive roots of $sl_{N+1}$ by
$\alpha_1,\dots,\alpha_N$.

\medskip
In our first example (Section~\ref{S41}), we consider arbitrary integral dominant weights
$\Lambda(1)$, $\Lambda(0)$ and assume $\vv_k$ to be a Bethe vector    
of the weight $\Lambda(1)+\Lambda(0)-k\alpha_1$.

\medskip
In another example (Section~\ref{S42}), we restrict 
$L_{\Lambda(0)}$ to be a symmetric power
of the standard $sl_{N+1}$-representation,  and
assume  $\vv_{k,1,1}$ to be a Bethe vector in $L$  of the weight
$\Lambda(1)+\Lambda(0)-k\alpha_1-\alpha_2-\alpha_3$.

\bigskip\noindent {\bf Theorem.}\ \  {\it Bethe vectors $\vv_k$ and
$\vv_{k,1,1}$ are non-trivial.}

\bigskip
For $N=1$,  any Bethe vector is  of the form $\vv_k$,
therefore the example~\ref{S41} implies  that  for $N=1$ and $n=2$ 
the Bethe vectors never vanish. Moreover, this example admits an
immediate generalization to $n>2$, see Theorem~\ref{T215} in 
Section~\ref{S415}. As a corollary we obtain
that if $L_{\Lambda(0)}=m\lambda_1$ and
$\Lambda(k_1,k_2)=\Lambda(1)+\Lambda(0)-k_1\alpha_1-k_2\alpha_2$ is
the highest weight of an irreducible component of
$L$, then Bethe vectors of the weight $\Lambda(k_1,k_2)$ do not vanish.
In particular, if $L_{\Lambda(1)}$ and $L_{\Lambda(0)}$ are 
$sl_3$-representations, then all Bethe vectors in $L$ do not vanish. 
 
\medskip
In the examples~\ref{S41}, \ref{S42} the subspace of singular vectors is
one-dimensional, therefore a Bethe vector, if it exists, gives an
eigenbasis of the Gaudin hamiltonians in the corresponding weight
subspace. A way to solve the Bethe equations in the example~\ref{S42}
is explained in~\cite{Sc}.

\bigskip
The key ingredient of our proof is to write the Bethe equations 
and projections of vector $\vv(\T)$ in terms of symmetric functions 
in $\T$, see Section~\ref{S51} and Section~\ref{S4}. 
Our calculations are based on funny relations between 
symmetric rational functions which generalize the ``Jacobi identity"
$$
\frac{1}{(x-y)(x-z)} + \frac{1}{(y-x)(y-z)} = \frac{-1}{(z-x)(z-y)}\,,
$$

\medskip\noindent
see Theorems~\ref{T1} and Corollary~ \ref{C1} in Section~\ref{S3}.

\bigskip\noindent{\bf Plan of the note.}\ \
Section~\ref{S3} is devoted to the ``Jacobi-like'' identities;
Section~\ref{S2} contains a description of the Bethe equations and
Bethe vectors;
in Section~\ref{S4} we verify non-triviality of  Bethe  vectors
in the examples.

\bigskip\noindent
{\bf Acknowledgments.}\ \
This work has been done in April--June 2004, when the second author visited 
the Mathematical Department of the Ohio State University.
It is a pleasure to thank this institution for hospitality and
excellent working  conditions.
We are also grateful to E.~Frenkel  and  V.~Lin for useful discussions,
and to A.Varchenko for his criticism of the first version of this note.
Finally, we would like to thank the referee for valuable remarks.

\bigskip
\section{Identities}\label{S3}

For  a function $g(t_1\dots t_k)$,  define {\it  symmetrization} as follows,
$$\s_k[g] :=  \sum_{\pi\in S_k}
    g\bigl(\pi(t_1),\dots, \pi(t_k)\bigr)\,,
$$
here the sum runs over the group $S_k$ of all permutations $\pi$ of variables
$t_1\dots t_k$.

\begin{thm}\label{T1} For any fixed $s_1$, $s_2$ and $s$, we have
$$
\begin{array}{l}\vspace{15pt}\displaystyle
{\rm(I)}\qquad
\s_k\left[\frac{1}{(s_1-t_1)(t_1-t_2)\dots(t_{k-1}-t_k)(t_k-s_2)}\right] = \\
\vspace{15pt}\displaystyle\ \qquad \ \qquad\ \qquad\ \qquad \
\qquad\ \qquad =\frac{(-1)^k\cdot(s_1-s_2)^{k-1}}{(s_1-t_1)
        \dots(s_1-t_k)\cdot(s_2-t_1)\dots(s_2-t_k)}\ ,\\
\vspace{15pt}\displaystyle {\rm(II)}\qquad
\s_k\left[\frac{1}{(s-t_1)(t_1-t_2)\dots(t_{k-1}-t_k) t_k}\right] =
  \frac{s^{k-1}}{(s-t_1)\dots(s-t_k)\cdot t_1\dots t_k}\ ,\\   \vspace{15pt}\displaystyle {\rm(III)}\qquad
\s_k\left[\frac{1}{(t_1-t_2)(t_2-t_3)\dots(t_{k-1}-t_k)(t_k-s)}\right] =
   \frac{(-1)^k}{(s-t_1)\dots(s-t_k)}\ .
\end{array}
$$
\end{thm}

\bigskip
{\bf Proof.}\ \ We prove the first and the third identities by
induction in $k$. The second identity can be obtained from the
first one by substitution  $s_1=s$ and $s_2=0$.

\medskip\noindent{\it The first identity}\ \ for $k=1$  becomes
$$
{\s}_1\left[\frac{1}{(s_1-t_1)\cdot(t_1-s_2)}\right] = \frac{1}{(s_1-t_1)
\cdot(t_1-s_2)} = \frac{-1}{(s_1-t_1)\cdot(s_2-t_1)}\,,
$$
and is true. Suppose that the identity (I) holds for $k-1$ and
prove it for $k$. Consider the subgroup $S_{k-1}\subset S_k$ of
permutations acting on the first $k-1$ variables. Every summand in
the symmetrization of our fraction has a form
$$
\frac{1}{(s_1-t_{i_1})(t_{i_1}-t_{i_2})\dots(t_{i_{k-1}}-t_{i_k})
      (t_{i_k}-s_2)}\ .
$$
Combine together all the summands with a fixed value of $i_k$, say
$i_k=j$, and factor out the last multiplier $1/(t_j-s_2)$. Then we
can write
$$
\begin{array}{l}\displaystyle
  \s_k\left[\frac{1}{(s_1-t_1)(t_1-t_2)\dots(t_{k-1}-t_k)(t_k-s_2)}\right]
      \quad = \vspace{10pt}\\   \displaystyle\hspace{4cm}
=\quad \sum_{j=1}^k \s_{k-1}\left[\frac{1}{(s_1-t_{i_1})(t_{i_1}-t_{i_2})
            \dots(t_{i_{k-1}}-t_j)}\right]\cdot\frac{1}{t_j-s_2}\ ,
\end{array}
$$
here the values of $i_1,\dots,i_{k-1}$ are different form $j$ and
the group $S_{k-1}$ acts by the permutations which keep $t_j$.
By the induction hypothesis this is
$$
\sum_{j=1}^k \frac{(-1)^{k-1} \cdot (s_1-t_j)^{k-2}}
{(s_1-t_1)\dots\wh{(s_1-t_j})\dots(s_1-t_k)\cdot
(t_j-t_1)\dots\wh{(t_j-t_j})\dots(t_j-t_k)}
   \cdot \frac{1}{t_j-s_2}\,,
$$
where the ``hat" means that the corresponding factor is omitted.
Multiplying this expression with
$(s_1-t_1)\dots(s_1-t_k)\cdot(s_2-t_1)\dots(s_2-t_k)$
we get
$$
(-1)^k\sum_{j=1}^k \left( (s_1-t_j)^{k-1}\cdot
   \frac{(s_2-t_1)\dots\wh{(s_2-t_j})\dots(s_2-t_k)}
        {(t_j-t_1)\dots\wh{(t_j-t_j})\dots(t_j-t_k)}\right)\ .
$$
This is nothing but the Lagrange interpolation formula for a
polynomial of degree $k-1$ in a variable $s_2$ that takes the
value $(-1)^k(s_1-t_j)^{k-1}$ at the point $s_2=t_j$ for every
$j=1,\dots,k$. Therefore it is equal to $(-1)^k(s_1-s_2)^{k-1}$.

\medskip\noindent{\it The third identity}\ \ is obvious for $k=1$,
$$
\s_1\left[\frac{1}{(t_1-s)}\right] = \frac{1}{(t_1-s)} =
\frac{-1}{(s-t_1)}\ .
$$
Suppose that the identity (II) holds for $k-1$ and prove it for
$k$. As before, combining together all the summands of the left
hand side with a fixed variable $t_j$ at the last factor of the
denominator we get
$$\begin{array}{l}\displaystyle
  \s_k\left[\frac{1}{(t_1-t_2)\dots(t_{k-1}-t_k)(t_k-s)}\right]
      \quad = \vspace{10pt}\\   \displaystyle\hspace{4cm}
=\quad  \sum_{j=1}^k \s_{k-1} \left[\frac{1}
{(s_1-t_{i_1})(t_{i_1}-t_{i_2})\dots(t_{i_{k-1}}-t_j)} \right]
\cdot \frac{1}{t_j-s}\ ,
\end{array}
$$
where the values of $i_1,\dots,i_{k-1}$ are different form $j$ and the group
$S_{k-1}$ acts by the permutations that keep $t_j$.
By the induction hypothesis
this is equal to
$$
\sum_{j=1}^k \frac{(-1)^{k-1}}{(t_j-t_{i_1}) (t_j-t_{i_2})\dots
    (t_j-t_{i_{k-1}})}\cdot \frac{1}{t_j-s}\,.
$$
Multiplying this expression with $(s-t_1)\dots(s-t_k)$
we get
$$
(-1)^k\sum_{j=1}^k\frac{(s-t_{i_1}) (s-t_{i_2})\dots(s-t_{i_{k-1}})}
  {(t_j-t_{i_1}) (t_j-t_{i_2}) \dots(t_j-t_{i_{k-1}})}\ ,
$$
where the indices $i_1,\dots,i_{k-1}$ in every summand are the integers
between $1$ and $k$ different form $j$. Recognizing in the last expression
the Lagrange interpolation formula we conclude that this is exactly $(-1)^k$.
\hfill$\square$

\bigskip
It is convenient to write identities on functions which are
symmetric with respect to variables $t_1,\dots, t_k$ in terms of
the elementary symmetry functions.

\bigskip\noindent{\bf Notation.}
$$
T(x)=(x-t_1)\,\dots\,(x-t_k)=x^k-\tau_1x^{k-1}+\, \dots\,+(-1)^k\tau_k\,,
$$
that is $\tau_i$ is the $i$-th elementary symmetric function in
$t_1,\dots, t_k$ for $1\leq i\leq k$; we set $\tau_0=1$.

\medskip
With this notation, the identities of Theorem~\ref{T1} take the
form
$$
\begin{array}{l}\vspace{15pt}\displaystyle
{\rm(I')}\qquad
\s_k\left[\frac{1}{(s_1-t_1)(t_1-t_2)\dots(t_{k-1}-t_k) (t_k-s_2)}\right] =
         \frac{(-1)^k\cdot(s_1-s_2)^{k-1}}{T(s_1)\cdot T(s_2)}\ ,\\
\vspace{15pt}\displaystyle {\rm(II')}\qquad
\s_k\left[\frac{1}{(s-t_1)(t_1-t_2)\dots(t_{k-1}-t_k) t_k}\right] =
         \frac{s^{k-1}}{T(s)\ \tau_k}\ ,\\
\vspace{15pt}\displaystyle {\rm(III')}\qquad
\s_k\left[\frac{1}{(t_1-t_2) (t_2-t_3)\dots(t_{k-1}-t_k)(t_k-s)}\right] =
         \frac{(-1)^k}{T(s)}\ .
\end{array}
$$

\begin{cor}\label{C1} We have
$$\begin{array}{l}\vspace{15pt}\displaystyle
{\rm(IV)}\qquad
\s_k\left[\frac 1{(t_1-t_2)\dots (t_{i-2}-t_{i-1})(t_{i-1}-s)
    (s-t_i)(t_i-t_{i+1})\dots (t_{k-1}-t_k)t_k}\right] = \\
\vspace{15pt}\displaystyle\ \hspace{8cm}
 = \frac{(-1)^{i-1}\ s^{k-i}\ \tau_{i-1}}{T(s)\ \tau_k}\,,
  \qquad \mbox{for\ } \quad 1\leq i\leq k\,.
\end{array}
$$
\end{cor}

\bigskip
{\bf Proof.}\ \ The formula (IV) for $i=1$ is exactly the identity
${\rm(II')}$. For $2\leq i\leq k$, first let us take the sum over
the subgroup $S_{i-1}\times S_{k+1-i}\subset S_k$, i.e. combine
together the summands corresponding to permutations of the first
$i-1$  and of the last $k+1-i$ variables $t_j$'s. Applying the
identities (III) and (II) we get
$$
\frac{(-1)^{i-1}}{(s-t_1)\dots(s-t_{i-1})}\times
\frac{s^{k-i}}{(s-t_i)\dots(s-t_k)\cdot t_i\dots t_k}\,.
$$
Now we take the sum over the cosets of the subgroup
$S_{i-1}\times S_{k+1-i}\subset S_k$ and collect look like terms.
Every denominator is the same,\quad
$(s-t_1)\dots(s-t_k)\cdot t_1\dots t_k = T(s)\cdot\tau_k$,\quad
whereas the numerators contain all possible products of $i-1$
of variables $t_j$'s. \hfill $\square$

\bigskip\noindent{\bf Remarks.}

{\bf 1.} Let us consider  $t_1,\dots, t_k$, $s_1$ as fixed numbers,
and $s$, $s_2$ as variables. Then the left-hand side
of every identity is nothing but a partial fraction decomposition
of the function from the right-hand side. This interpretation,
indicated by V.~Lin, leads to another proof of the identities.

{\bf 2.} As A.~Varchenko pointed out,  our identity (III) for $s=0$ follows
from the coincidence of the forms $\Omega^{sl_2}$ and $\wt{\Omega}^{sl_2}$ 
from~\cite[page~2] {RStV}.  Notice that for arbitrary
$s$  the identity (III) can be obtained from this particular case by
substitution $\T\mapsto \T-s$.
Similarly, the substitution $s\mapsto s_1-s_2$, $\T\mapsto \T-s_2$
transforms the identity (II) into the identity (I).

{\bf 3.} A remark of the referee  is that  the identity (III) could be 
deduced from the identity (I).  Indeed if  we consider  (I)  as
a function of a complex variable $s_1$ and take the residues of both sides
at infinity, then we  get  (III).

\section{Bethe vectors}\label{S2}
Here we recall the constructions for the tensor product
of $n=2$ modules corresponding to points $z_0=0$ and $z_1=1$.
For $n>2$ (and for any simple Lie algebra), see \cite{FeFRe, ReV}.

\subsection{Subspace of singular vectors in $L$}\label{S20} \ \
Denote by $\{e_i,f_i,h_i\}_{i=1}^{N}$
the standard Chevalley generators of $sl_{N+1}(\C)$,
$$
[h_i,e_i]=2e_i\,,\ [h_i,f_i]=-2f_i\,,\ [e_i,f_i]=h_i\,;\qquad
[h_i,h_j]=0\,,\ [e_i,f_j]=0\,  {\rm \ if \ } \ i\neq j\,.
$$
Let $\h$ be the Cartan subalgebra and $\h^*$ its dual,
$$
\h^*\,=\,\C\{\lambda_1,\dots,\lambda_{N+1}\}\,/\,
\left( \lambda_1+\dots+\lambda_{N+1}=0 \right)\,,
$$
with the standard bilinear form $(\cdot,\cdot)$.
The simple positive roots are  $\alpha_i=\lambda_i-\lambda_{i+1}$,
$1\leq i\leq N$,
$$
(\alpha_i,\alpha_i)=2\,;\quad
(\alpha_i,\alpha_j)=0\,, \mbox{\ if\ } |i-j|>1\,;\quad \mbox{and}\quad
(\alpha_i,\alpha_j)=-1\,, \mbox{\ if\ } |i-j|=1\,.
$$

\medskip
Let $\Lambda(1)$ and $\Lambda(0)$ be integral dominant weights,
and $\kk=(k_1,\dots, k_N)$ be a vector with nonnegative integer
coordinates such that
$$\Lambda(\kk):=\Lambda(1) + \Lambda(0)-k_1\alpha_1-\dots-k_N\alpha_N
$$
is an integral dominant weight as well. Denote by
$$\Sing_{\kk}L:=\{\vv\in L\,
\vert\, h_i\vv=(\Lambda(\kk),\alpha_i)\vv,\ e_i\vv=0, \ i=1,\dots,N \}
$$
the subspace of singular vectors in
$L=L_{\Lambda(1)}\,\otimes\,L_{\Lambda(0)}$  of  weight $\Lambda(\kk)$.

\subsection{Bethe system associated with $\Sing_{\kk}L$}\label{S21}
For every $i=1,\dots, N$ introduce a set of $k_i$ auxiliary variables
associated with the root $\alpha_i$,
$$
t(i):= \left(\, t_1(i)\,,\  \dots\,,\  t_{k_i}(i)\right)\,,
$$
and  write $\T:=\left(\, t(1)\,,\ \dots\,, t(N)\,\right)$.

\medskip
The {\it Bethe system} is the following system of equations on variables
$t_l(i)$,
$$\sum_{s\neq l}\frac 2{t_l(i)-t_s(i)}-
  \sum_{s=1}^{k_{i-1}}\frac 1{t_l(i)-t_s(i-1)}-
  \sum_{s=1}^{k_{i+1}}\frac 1{t_l(i)-t_s(i+1)}-
  \frac {(\Lambda(0),\alpha_i)}{t_l(i)} -
  \frac {(\Lambda(1),\alpha_i)}{t_l(i)-1}\,=\,0\,,
$$
here  $1\leq i\leq N,\ 1\leq l\leq k_i$.

\medskip
Every solution $\T^{(0)}$ to this system 
determines a {\it Bethe vector} $\vv(\T^{(0)})=\vv_{\kk}(\T^{(0)})\in \Sing_{\kk}L$. 
The function $\vv(\T)$  is described   in Section~\ref{S22}.

\bigskip
\subsection{Bethe equations in  terms of polynomials
$T_1(x),\dots, T_N(x)$.}\label{S51}
We use the notation introduced in Section~\ref{S3}.

\begin{prop}\label{P1}
Assume all the roots of $T(x)$ to be simple. Then
$$
\frac{T'(x)}{T(x)}=\sum_{j=1}^k\ \frac 1{x-t_j}\,;\quad
\frac{T''(t_i)}{T'(t_i)}=\sum_{j\neq i}\ \frac 2{t_i-t_j}\,.
$$
\end{prop}

\bigskip{\bf Proof.}\ \ The first equation is just
the logarithmic derivative of $T$.

We have
$$
T'(x)\,=\,\left(\ \sum_{j=1}^k\ \frac 1{x-t_j}\ \right)\cdot T(x)\,.
$$
Derivation of this equation gives
$$
T''(x)\,=\,\left(\ \sum_{j=1}^k\ \frac 1{x-t_j}\ \right)'\cdot T(x)+
\left(\ \sum_{j=1}^k\ \frac 1{x-t_j}\ \right)^2\cdot T(x)\,.
$$
Therefore
$$
\frac{T''(x)}{T(x)}\,=\, -\,\sum_{j=1}^k\ \frac 1{(x-t_j)^2}\,+\,
\left(\ \sum_{j=1}^k\ \frac 1{x-t_j}\ \right)^2\,=\,
2\sum_{1\leq j<l\leq k}\ \frac 1{(x-t_j)(x-t_l)}\,.
$$
We have
$$\frac{T''(x)}{T'(x)}\,=\,
  \frac{\displaystyle 2\sum_{1\leq j<l\leq k}\ \frac 1{(x-t_j)(x-t_l)}}
       {\displaystyle \sum_{j=1}^k\ \frac 1{x-t_j}}\,  =
  \frac{\displaystyle 2\sum_{1\leq j<l\leq k}(x-t_1)\dots
                    \wh{(x-t_j})\dots \wh{(x-t_l})\dots(x-t_k)}
{\displaystyle \sum_{j=1}^k (x-t_1)\dots \wh{(x-t_j})\dots (x-t_k)}\,.
$$
Substitution $x=t_i$ gives
$$\frac{T''(t_i)}{T'(t_i)}\,=\,
  \frac{\displaystyle 2\sum_{j\neq i}(t_i-t_1)\dots
        \wh{(t_i-t_j})\dots \wh{(t_i-t_i})\dots(t_i-t_k)}
       {(t_i-t_1)\dots \wh{(t_i-t_i})\dots (t_i-t_k)}\,,
$$
and the division finishes the proof. \hfill $\square$

\bigskip
Now we can re-write the Bethe system
in terms of polynomials $T_1(x),\dots, T_N(x)$, where
$$
T_i(x)=\left(x-t_1(i)\right)\,\dots\, \left(x-t_{k_i}(i)\right)\,.
$$
We have
$$
 \frac {T_i''\left(t_l(i)\right)}{T_i'\left(t_l(i)\right)}
 -\frac{T_{i-1}'\left(t_l(i)\right)}{T_{i-1}\left(t_l(i)\right)}
 -\frac{T_{i+1}'\left(t_l(i)\right)}{T_{i+1}\left(t_l(i)\right)}
 -\frac {(\Lambda(0),\alpha_i)}{t_l(i)} -
\frac {(\Lambda(1),\alpha_i)}{t_l(i)-1}\,=\,0\,,
$$
for   $1\leq i\leq N,\ 1\leq l\leq k_i$.

\bigskip
\subsection{Function $\vv(\T)$}\label{S22}
The function $\vv(\T)$ has been obtained in~\cite[Sections~6,7]{SV}
in general setting (see also~\cite{MV1}, where it is called 
{\it the universal weight function}). Below we rewrite this function for 
the weight $\Lambda(\kk)$ in the tensor product of two modules
$L=L_{\Lambda(1)}\otimes L_{\Lambda(0)}$ corresponding to points
$z_1=1$ and $z_0=0$. Generic case can be found in~\cite{FeFRe, ReV}.

\medskip
To begin with, we construct the vectors that  generate the subspace
$L_{\Lambda(\kk)}\subset L$ of weight $\Lambda(\kk)$.
In general, their number is greater than the dimension of that subspace,
so they are linearly dependent.

\medskip
Consider all  pairs of words $(\bF_1,\bF_0)$ in letters $f_1,\dots, f_N$
subject
to the condition that the total number of occurrences of letter $f_i$ in 
both words is precisely $k_i$. Our vectors will be labeled by these pairs.
Namely, we may think about words $\bF_1$ and $\bF_0$ as
elements of the universal enveloping algebra of $sl_{N+1}$ that naturally
act on the spaces $L_{\Lambda(1)}$ and $L_{\Lambda(0)}$, respectively.
Fix the highest weight vectors\
$\vv_1\in L_{\Lambda(1)}\,,\ \vv_0\in L_{\Lambda(0)}$.
Then the vector
$$
 \w_{_{(\bF_1,\bF_0)}} := \bF_1\vv_1\,\otimes\,\bF_0\vv_0
$$
has weight $\Lambda(\kk)$, and all such vectors generate the weight space
$L_{\Lambda(\kk)}$.

\medskip
Now we define $\vv_{\kk}(\T)$ as a linear combination
$$
\vv_{\kk}(\T)\ :=\ \sum_{(\bF_1,\bF_0)}\
      \omega_{_{(\bF_1,\bF_0)}}(\T)\,\w_{_{(\bF_1,\bF_0)}}\,,
$$
where $\omega_{_{(\bF_1,\bF_0)}}(\T)$ are certain rational functions.
We will construct these functions in two steps described below.
Write
$$\bF_1 = f_{i_1}\dots f_{i_{s_1}}\,,\qquad
  \bF_0 = f_{j_1}\dots f_{j_{s_0}}\,, \qquad
  \bF_1\,\bF_0 = f_{i_1}\dots f_{i_{s_1}}\,f_{j_1}\dots f_{j_{s_0}}\,.
$$
The length of the word $\bF_1\bF_0$ equals $s_1+s_0=k_1+\dots+k_N$.

\medskip
{\it The first step} is to translate $(\bF_1,\bF_0)$ into
a rational function
$g_{_{(\bF_1,\bF_0)}}(\T)$ of $\T$.
For every  $i=1,\dots, N$, we replace the
first occurrence (from left to right) of $f_i$ in the word
$\bF_1\bF_0$ by the variable $t_1(i)$; the second occurrence 
by the variable $t_2(i)$; and so on up to the last, $k_i$-th,
occurrence, where $f_i$ will be replaced by $t_{k_i}(i)$. We will get
a pair of words in $\T$.  We augment these two words by
1 and 0, according to the values of $z_1$ and $z_0$,
and thus get the row,
$$t_{a_1}(i_1)\, t_{a_2}(i_2)\,\dots\, t_{a_{s_1}}(i_{s_1})\,1\,,\qquad
  t_{b_1}(j_1)\, t_{b_2}(j_2)\,\dots\, t_{b_{s_0}}(j_{s_0})\,0\,,
$$
in which every variable $t_l(i)$ from $\T$ appears precisely once.
This row defines the fraction
$$\begin{array}{l} \displaystyle
g_{_{(\bF_1,\bF_0)}}(\T) :=
\frac1{\bigl(t_{a_1}(i_1)-t_{a_2}(i_2)\bigr)
       \bigl(t_{a_2}(i_2)-t_{a_3}(i_3)\bigr) \cdots
       \bigl(t_{a_{s_1-1}}(i_{s_1})-t_{a_{s_1}}(i_{s_1})\bigr)
       \bigl(t_{a_{s_1}}(i_{s_1})-1\bigr)}
    \vspace{15pt} \\ \hspace{2.8cm}\displaystyle \times\
 \frac1{\bigl(t_{b_1}(j_1)-t_{b_2}(j_2)\bigr)
        \bigl(t_{b_2}(j_2)-t_{b_3}(j_3)\bigr) \cdots
       \bigl(t_{b_{s_0-1}}(j_{s_0})-t_{b_{s_0}}(j_{s_0})\bigr)
       t_{b_{s_0}}(j_{s_0})}\,.
\end{array}
$$

\medskip
{\it The second step} is the symmetrization of $g_{_{(\bF_1,\bF_0)}}(\T)$.
Let $S_{\kk}$ denote the group of permutations of variables
$$\T=\left(t_1(1),\dots, t_{k_1}(1),\ t_1(2),\dots, t_{k_2}(2),
    \ \dots\ ,\ t_1(N),\dots, t_{k_N}(N)\right)
$$
that permute variables $t_1(i),\dots,t_{k_i}(i)$ within their own,
$i$-th, set,
for every $i=1,\dots,N.$
Thus $S_{\kk}$ is isomorphic to the direct product
$S_{k_1}\times S_{k_2}\times\dots\times S_{k_N}$ of permutation groups.

For a function $g(\T)$ define
the symmetrization by the formula
$$\s_{\kk}[g] :=  \sum_{\pi\in S_{\kk}}
g\Bigl(\pi\bigl(t_1(1)\bigr),\dots, \pi\bigl(t_{k_N}(N)\bigr)\Bigr)\,.
$$

Finally we set
$$\omega_{_{(\bF_1,\bF_0)}}(\T) :=
\s_{\kk}\left[g_{_{(\bF_1,\bF_0)}}(\T)\right]\,.
$$

Notice that the universal weight function $\vv_{\kk}(\T)$ is defined
for any, not necessarily dominant, weight $\Lambda(\kk)$ presented in
$L=L_{\Lambda(1)}\otimes L_{\Lambda(0)}$.
However in the Bethe Ansatz it is used  only when $\Lambda(\kk)$ is
the highest weight of an irreducible component of $L$.

\section{Checking the non-triviality of Bethe vectors in
         examples}\label{S4}

\subsection{Example $\Lambda(1)+\Lambda(0)-k\alpha_1$}\label{S41}\ \
We assume $\Lambda(1)$ and $\Lambda(0)$ to be integral dominant
weights and $k$ be an integer such that 
$$
\Lambda(k,0)\,:=\,\Lambda(1)+\Lambda(0)-k\alpha_1
$$
is the highest weight of an irreducible component of
$L=L_{\Lambda(1)}\otimes L_{\Lambda(0)}$. 
In this case Steinberg's formula
implies that $f_1^k\,\vv_0 \neq \0$ \cite[Exercise 24.12]{Hu}.
We will show that the universal 
weight function $\vv_k(\T):=\vv_{(k,0,\dots,0)}(\T)$
never vanish. We have
$$
\kk=(k,0,\dots,0)\,,\quad \T=\left(t(1)\right)\,,\quad \bF_0=f_1^k\,,\quad
g_{_{(\emptyset,\bF_0)}}(\T)=
     \frac 1{(t_1-t_2)\dots (t_{k-1}-t_k)t_k}\,.
$$
Simplifying the notation, write
$$ t:=(t_1,\dots,t_k)\,,\qquad T(t)=\prod_{i=1}^k(x-t_i)\,,\qquad
  \omega_{k,0}(t) := \omega_{_{(\emptyset,\bF_0)}}(\T)\,.
$$

\begin{thm}\label{T2}
The projection of the vector $\vv_k(\T)$
to the subspace of $L_{\Lambda(k,0)}$ spanned by
$\vv_1\,\otimes\,f_1^k\,\vv_0$ is
a non-zero vector.
\end{thm}

\bigskip
{\bf Proof.}\ \ Notice that the domain of the function
$\vv_k(\T)$ is given by the inequalities,
$$
t_i\neq t_j\,,\quad t_i\neq 0\,,\quad 1\leq i\neq j\leq k\,.
$$
The considered projection of $\vv_k(\T)$  has the form
$\omega_{k,0}(t)\vv_1\,\otimes\,f_1^k\,\vv_0$, where
$$
\omega_{k,0}(t)\,=\,\s_k \left[\frac 1{(t_1-t_2)\dots
(t_{k-1}-t_k)t_k} \right]\, .
$$
The identity (III) with $s=0$ gives
$$
\omega_{k,0}(t)\,=\,\frac{(-1)^k}{T(0)}\,=\,\frac 1{t_1\, t_2\dots
t_k}\,,
$$
and this fraction never vanishes.  \hfill$\square$

\

\subsection{Generalization to arbitrary $n$}\label{S415}\ \
Theorem~\ref{T2} has the following generalization to the universal
weight function $\vv(\T)$ corresponding to the weight
$$
\Lambda(k_1,\dots,k_m)=\sum_{i=1}^n \Lambda(i)\, -
\,\sum_{i=1}^m k_i\alpha_i
$$
in the tensor product
$$
L=L_{\Lambda(1)}\otimes\dots\otimes L_{\Lambda(n)}
$$
of $n$ highest weight  $sl_{N+1}$-representations marked by
distinct complex numbers $z_1,\dots,z_n$. 

\begin{thm}\label{T215}
Assume  that $m\leq \min(n,N)$ and
$$
k_i \leq (\Lambda(i),\alpha_i)\,, \quad i=1,\dots,m\,.
$$
Then the universal weight function $\vv(\T)$ corresponding to the weight 
$\Lambda (k_1,\dots,k_m)$ does not vanish.
\end{thm}

\bigskip
{\bf Proof.}\ \ Fix  highest weight 
vectors $\vv_i\in L_{\Lambda(i)}\,,$\ \  $i=1,\dots,n$. 
According to our assumptions, we have
$$
f^{k_i}_i\vv_i\neq \0\,, \quad i=1,\dots,m.
$$
Consider the projection of $\vv(\T)$ to the
one-dimensional subspace of $L$ spanned by
$$
f^{k_1}_1\vv_1\otimes \dots \otimes f^{k_m}_m\vv_m \otimes
\vv_{m+1} \otimes  \dots \otimes \vv_n\,.
$$
Applying the identity (III), one gets that the corresponding coefficient
is equal to
$$\frac{(-1)^{k_1+\dots +k_m}}{T_1(z_1)\dots T_m(z_m)}\ ,
$$
where polynomials $T_i(x)$ are as in Section~\ref{S51}, and hence
does not vanish.   \hfill$\square$

\medskip
In particular, if the Bethe vector of the weight $\Lambda (k_1,\dots,k_m)$
exists, then it is a non-zero vector.  

\medskip
Returning to $n=2$, in the case when one of  two modules
is a symmetric power of the standard representation, we arrive at the 
following result.

\begin{cor}
If $L_{\Lambda(0)}=m\lambda_1$ and
$\Lambda(k_1,k_2)=\Lambda(1)+\Lambda(0)-k_1\alpha_1-k_2\alpha_2$ is
the highest weight of an irreducible component of
$L=L_{\Lambda(1)}\otimes L_{\Lambda(0)}$, then the universal weight function
$\vv(\T)$ corresponding to the weight $\Lambda(k_1,k_2)$ does not vanish.
In particular, if $L_{\Lambda(1)}$ and $L_{\Lambda(0)}$ are 
$sl_3$-representations, then all Bethe vectors in $L$ do not vanish. 
\end{cor}

\bigskip
{\bf Proof.}\ \
Elementary considerations with the Pieri formula
\cite[Proposition 15.25]{FuH} show that the conditions
$$
k_1\leq (\Lambda(0),\alpha_1)\,, \quad k_2\leq (\Lambda(1),\alpha_2)
$$
are always fulfilled. In the $sl_3$ case all highest
weights are clearly of the form $\Lambda(k_1,k_2)$.\hfill$\square$

\subsection{Example $\Lambda(1)+\Lambda(0)-k\alpha_1-\alpha_2-\alpha_3$}
\label{S42}\ \
We assume that $\Lambda(0)=m\lambda_1$,\ $N \geq 3$, and $k\geq 1$ is an integer such that
$$\Lambda(k,1,1)\,:=\,
\Lambda(1)+\Lambda(0)-k\alpha_1-\alpha_2-\alpha_3
$$
is the highest weight of an irreducible component of
$L=L_{\Lambda(1)}\otimes L_{\Lambda(0)}$.
As before, the Pieri formula \cite[Proposition 15.25]{FuH}
implies that $f_2f_1^k\vv_0 \neq \0$ and 
$f_3\vv_1 \neq \0$ for any fixed highest weight
vectors $\vv_0\in L_{\Lambda(0)}$ and $\vv_1\in L_{\Lambda(1)}$.

\medskip
 The module $L_{\Lambda(0)}$ is the $m$-th symmetric
power of the standard $sl_{N+1}$-representation. Take
$\vv_0=\epsilon_1^m$, where $\{\epsilon_i\}$ is a basis in the
standard representation,
$$
f_i\epsilon_i=\epsilon_{i+1}\,,\quad f_i\epsilon_j=\0\,,\quad  i\neq j\,.
$$
The subspace of weight $\Lambda(0)-k\alpha_1-\alpha_2$ in $L_{\Lambda(0)}$
is one-dimensional and generated by the vector
$\epsilon_1^{m-k}\,\epsilon_2^{k-1}\,\epsilon_3$.

There are three sets of auxiliary variables. We write
$$t(1)=(t_1,\dots,t_k)\,,\ t(2)=s\,,\ t(3)=r\,,\
\T=(t,s,r)=(t_1,\dots,t_k,s,r)\,, \  T(x)=\prod_{i=1}^k\,(x-t_i)\,.
$$
The word $\bF_0$ can be written as
$\bF_0=f_1^{i-1}f_2f_1^{k+1-i}$ for $i=1,\dots,k+1$.
Notice that $f_1^k f_2 \vv_0 = \0$ for our choice of
$\Lambda(0)$, therefore we  assume that $i$ varies from $1$ to $k$
and set
$$\omega_i(t,s,r) :=
  \omega_{_{\left(f_3,f_1^{i-1}f_2f_1^{k+1-i}\right)}}(\T)\,,
\quad i=1,\dots,k.
$$

\begin{thm}\label{T3}
If $T'(s)\neq 0$, then the projection of the vector $\vv_{k,1,1}(t,s,r)$ to the subspace
$L_{\Lambda(k,1,1)}$ spanned by
$$
f_3\vv_1\,\otimes \left(\epsilon_1^{m-k}\,\epsilon_2^{k-1}\,
\epsilon_3\right)
$$
is a non-zero vector.
\end{thm}

\bigskip
{\bf Proof.}\ \
Notice that the domain of the function $\vv_{k,1,1}(t,s,r)$ is given by
the inequalities,
$$
t_i\neq t_j\,,\quad, t_i\neq s\,,\quad, t_i\neq r\,,\quad, r\neq s\,,
\quad  t_i,s,r\neq 0,1\,, \quad 1\leq i\neq j\leq k.
$$
The projection of  vector $\vv_{k,1,1}(t,s,r)$ to the chosen subspace
has the form
$$
\sum_{i=1}^k\, \omega_i(t,s,r)\w_i\,,
$$
where $\w_i=f_3\vv_1\,\otimes\,f_1^{i-1}f_2\,f_1^{k+1-i}\vv_0$.
Here $f_2$ stands at the $i$-th place from the left, and
$$
\omega_i(t,s,r)=\frac{1}{r-1}
\s_k\!\left[\frac 1{(t_1-t_2)\!\dots\!
(t_{i-2}-t_{i-1})(t_{i-1}-s)
(s-t_i)(t_i-t_{i+1})\!\dots\!(t_{k-1}-t_k)t_k}\right].
$$

\medskip
An easy calculation shows that
$$
\w_i\,=\,(k+1-i)\cdot m(m-1)\dots (m+1-k)\cdot
f_3\vv_1\,\otimes \left(\epsilon_1^{m-k}\,\epsilon_2^{k-1}\,\epsilon_3\right)\,,
$$
therefore the projection is
$$
\left(\,\sum_{i=1}^k\,(k+1-i)\cdot \omega_i\,(t,s,r)\right)
\cdot m(m-1)\dots (m+1-k)\cdot
f_3\vv_1\,\otimes \left(\epsilon_1^{m-k}\,\epsilon_2^{k-1}\,\epsilon_3\right)\,.
$$
It is convenient to use the notation introduced at the end of Section~\ref{S3}.
The identity (IV) of Corollary~\ref{C1} gives
$$
\omega_i(t,s,r)=
    \frac{(-1)^{i-1}\ s^{k-i}\ \tau_{i-1}}{(r-1)T(s)\ \tau_k}\,,
\quad 1\leq i\leq k\,.
$$
Therefore
$$
\sum_{i=1}^k\,(k+1-i)\cdot \omega_i\,(t,s,r)\,=\,
\frac{T'(s)}{(r-1)T(s)\tau_k}\,,
$$
and the statement of the theorem follows.  \hfill $\square$

\medskip
\begin{cor} If the Bethe vector $\vv_{k,1,1}$ exists, then it does not vanish.
\end{cor}

\bigskip
{\bf Proof.}\ \
We show that  $T'(s)$ can not vanish at a solution of the Bethe system.
The Bethe equation corresponding to the variable $s$ is as follows,
$$
\frac1{s-r} + \frac{T'(s)}{T(s)}+\frac{(\Lambda(1),\alpha_2)}{s-1}=0\,,
$$
whereas the one corresponding to $r$ has the form
$$
\frac1{r-s} + \frac{(\Lambda(1),\alpha_3)}{r-1}=0\,,
$$
see Section~\ref{S51}. Denote $(\Lambda(1),\alpha_2)=A$ and
$(\Lambda(1),\alpha_3)=B$. Assuming $T'(s)=0$ one gets  
the following linear system with respect to $s$ and $r$,
$$
-Ar+(A+1)s=1\,,\quad (B+1)r-Bs=1\,.
$$
The solution to this system is $r=s=1$ and contradicts to the
conditions $r\neq s\neq 1$.        \hfill $\square$

\newpage
\

\bigskip

\

\vskip1.cm

\parbox[t]{3in}{\it \text{Sergei~Chmutov}\\
The Ohio State University, Mansfield\\
1680 University Drive\\
Mansfield, OH 44906\\
USA\\
~\texttt{chmutov@math.ohio-state.edu}} \hspace{1.5cm}
\parbox[t]{2.5in}{\it \text{Inna~Scherbak}\\
School of Mathematical Sciences\\
Tel Aviv University\\
Ramat Aviv, Tel Aviv 69978\\
Israel\\
~\texttt{scherbak@post.tau.ac.il}}

\end{document}